\documentclass[a4paper, 10pt, conference]{ieeeconf}      

\IEEEoverridecommandlockouts                              

\usepackage{inputenc}

\usepackage{pgfplots}
\pgfplotsset{compat=1.18,width=\textwidth}

\usepackage{amsmath}
\usepackage{amsfonts}
\usepackage{amssymb}
\usepackage{nicefrac}
\usepackage{mathtools}

\usepackage{siunitx}
\usepackage{graphicx}
\usepackage[ruled,vlined]{algorithm2e}
\SetKw{Break}{break}

\usepackage{optidef}
\usepackage{comment}
\usepackage{csquotes}
\usepackage{float}

\usepackage{xcolor}

\usepackage{caption}
\usepackage{subcaption}
\usepackage{dirtytalk}
\usepackage{layouts}

\usepackage{booktabs}

\usepackage{svg}

\usepackage{mathrsfs}
\usetikzlibrary{arrows,angles,quotes,calc}

\usepackage{hyperref}
\hypersetup{
    colorlinks=true,
    linkcolor=red,
    filecolor=magenta,      
    urlcolor=blue,
    citecolor=purple,
    }

\usepackage{listings}
\usepackage[
    backend=biber,
    style=ieee,
    sorting=none]{biblatex}
\graphicspath{{figures/}}

\NewDocumentCommand{\setR}{}{\mathbb{R}}
\NewDocumentCommand{\setS}{}{\mathbb{S}}
\NewDocumentCommand{\setC}{}{\mathbb{C}}

\NewDocumentCommand{\norm}{m}{\left\lVert#1\right\rVert}
\NewDocumentCommand{\abs}{m}{\left\lvert#1\right\rvert}

\NewDocumentCommand{\brackets}{m}{\left[#1\right]}

\NewDocumentCommand{\parentheses}{m}{\left(#1\right)}
\NewDocumentCommand{\dv}{mm}{\frac{\mathrm{d}#1}{\mathrm{d}#2}}
\NewDocumentCommand{\bigO}{m}{\mathcal{O}\parentheses{#1}}

\RenewDocumentCommand{\refeq}{m}{(\ref{#1})}
\NewDocumentCommand{\citetodo}{O{}}{\texttt{[\textbackslash cite TODO #1]}}

\title{\LARGE \bf
Safe and wind-aware synchronous path planning for a fleet of fixed-wing constant speed aircraft
}

\author{Maël Feurgard$^{1}$, Gautier Hattenberger$^{1}$, Nicolas Durand$^{1}$ and Simon Lacroix$^{2}$
\thanks{$^{1}$ Fédération ENAC ISAE-SUPAERO ONERA, Université de Toulouse, France %
    {\tt\small {[name]}.{[surname]}@enac.fr}}%
\thanks{$^{2}$ Deceased, was with Laboratoire d'Analyse et d'Architecture des Systèmes, Université de Toulouse, CNRS, Toulouse, France.
    \url{https://sites.laas.fr/Simon_Lacroix/}}%
}

\begin{document}

\maketitle
\thispagestyle{empty}
\pagestyle{empty}

\begin{abstract}

    Path planning for multiple unmanned aerial vehicles is a difficult task, and even more for a fleet of fixed-wing aircraft.
    One specific case is the transition to, or between, formation flight patterns,
    which requires synchronous arrivals while ensuring minimal separation, and ideally maintaining cruise speed.
    We present a centralized method to solve this problem based on enumerating different Dubins paths.
    Given a travel time for the fleet, it builds a set of possible paths for each aircraft. Then,
    it checks in parallel separation between each path pair.
    This yields coefficients for an Integer Linear Programming problem determining if a fleet-wide conflict-free solution exists.
    This process is repeated for different travel times sampled with increasing resolution until
    the user-defined accuracy is met.
    The method is benchmarked with Monte-Carlo simulations considering up to 20 aircraft simultaneously,
    achieving an 95\% success rate for an average 8 seconds computation time.

\end{abstract}

\section{Introduction}


Unmanned Aerial Vehicles (UAVs) are getting more and more available at relatively low cost,
motivating both applications and research flying several autonomous aircraft simultaneously.
In particular, the deployment of several fixed-wing aircraft is relevant to long missions, requiring
significant payload or airspeed ; aerial cartography \cite{Salami2014}, airborne particles sensing \cite{Gonzalez2011,Techy2008}
and forest monitoring \cite{Casbeer2005,Salinas2023} are examples of such tasks.
They all require solving a common problem: cooperative path planning. 

Cooperative path planning consists in moving multiple agents in a shared space to achieve a common goal.
Several constraints can be considered, the most common being collision avoidance, synchronous arrival to destination,
maintain formation and feasible dynamics. This problem has been significantly studied for aerial vehicles \cite{PathPlanningReview,UAVSwarmReview,FormationControlReview},
with some methods dating back from the 1990s \cite{Pellazar1998}. 
We are particularly interested in collision-free 2D movements with synchronous arrivals, since such movements correspond to
transitions between planar formations which may be needed for take off, landing or during the mission. Different strategies exist
to achieve such transitions: discretize the airspace using a grid \cite{Pellazar1998,Shorakaei2014,AStarThenBezierPP},
a Voronoi decomposition \cite{McLain2000} or sets of lines \cite{Babel2018} first, then generate safe flight paths and finally smooth them;  
use leader-follower or swarm schemes to handle collision avoidance dynamically, and vary speed for synchronization 
\cite{Chen2021a,Wang2019a,Liu2022,Song2023,Madabushi2025}; or directly generate valid paths.

We focus on this last category as it does not rely on a discretization parameter and provide an explicit path to be followed for each agent.
This ensures that there are no deadlocks, but usually requires more computational power and may still fail to find a solution in some cases.
Paths can be generated along several models; some can consider the acceleration limitations of the aircraft, and will use clothoids \cite{Shanmugavel2007},
Bézier splines \cite{Yu_2017,DubinsToBSplineSync} or Pythagorean Hodograph splines \cite{Shanmugavel2007a,Cichella2014,Choe2016}; 
others will use the Dubins model \cite{Dubins1957},
as it simplifies trajectory generation into a geometry problem with circles and lines \cite{Shanmugavel2005,Chen2018} at the cost of
introducing acceleration discontinuities.

Among the full path generation methods, only a limited amount of aircraft are considered and performances are poorly studied. In \cite{DubinsToBSplineSync},
several Monte Carlo simulations are done, measuring accuracy and runtime of the planner, testing from 2 to 7 aircraft. But this article
does not manage collisions. The other articles show only sample cases, without timing computations and the number of aircraft does not exceed 5.
This is especially limited compared to the leader follower techniques, such as in \cite{Liu2022} where 3 groups of 7 aircraft were deployed in real conditions.


We present a new method based on \emph{Dubins paths} to solve the
fixed-wing fleet path planning problem with the constraints of constant airspeed, collision avoidance
and synchronous arrival in 2D. The main difference with \cite{Shanmugavel2005,Chen2018} is that instead
of looking for one specific combination of curves in a continuous domain satisfying both synchronization and collision avoidance,
we build sets of synchronous curves for each aircraft in which we look for a collision-free combination. This allows us to parallelize computations
and push the benchmark from 3 to 20 aircraft simultaneously. Furthermore, the method is able to take into account
a constant uniform wind.


\section{Problem formulation}
\label{sec:problem}

We consider the problem of path planning for a fleet of fixed-wing aircraft flying at constant air speed and altitude, using
simplified dynamics in an obstacle-free area.

Let the aircraft be numbered from 1 to $N$, and represented using a 2D position $[x_k,y_k]\in\setR^2$ and a 2D orientation 
$\theta_k\in\setS \coloneq (-\pi,\pi]$ with respect to the ground. We combine these values in a state vector $p_k = [x_k,y_k,\theta_k]$.
For the sake of simplicity, we assume all aircraft have the same characteristics, namely the
same cruise air speed $V$, minimal turn radius $\rho_\text{min}$, and minimal separation distance $\delta$.

This yields the following dynamics for each aircraft $k\in[1,N]$, the so-called \emph{Dubins vehicle} \cite{Dubins1957}:
\begin{equation}
    \label{eq:dubins}
    \dv{p_k}{t} = 
    \begin{pmatrix}
        \dot{x_k}\\
        \dot{y_k}\\
        \dot{\theta_k}
    \end{pmatrix} \coloneq \begin{pmatrix}
        V \cos \theta_k\\
        V \sin \theta_k\\
        u_k
    \end{pmatrix}
\end{equation}
The sole control parameter is the change in orientation $u_k\in\brackets{\nicefrac{-V}{\rho_\text{min}},\nicefrac{V}{\rho_\text{min}}}$. 
Note that the vehicle is always aligned with its speed vector. As such, for a differentiable trajectory \(\gamma\ \colon\ \setR \to \setR^2\),
we may abuse the notation $p_k(t) = \gamma(t)$ to state that $\gamma(t) = [x_k(t),y_k(t)]$ and 
$\gamma'(t)$ is aligned with $\brackets{\cos\theta_k(t) , \sin\theta_k(t)}$, 
i.e. the aircraft follows $\gamma$ both in position and orientation.


We require synchronicity: all aircraft reach their destination at the same time. This constraint
stems from formation flights, and can be adapted to sequencing for instance.
We develop this aspect in section \ref{sec:fleetpp}.

Let us formulate our optimization problem. Given for each aircraft $k$ an initial state $S_k$,
a final state $E_k$ both in $\setR^2 \times \setS$, and a set of possible twice differentiable paths 
$\Gamma\subset\mathcal{D}^2(\setR^+,\setR^2)$, we seek to minimize
the time taken by the fleet to reach its final configuration while satisfying the constraints. Formally:
\begin{equation}
    \label{eq:optim}
    \min_{
       \begin{matrix}
           \scriptstyle\tau \in\, \setR^+\\
           \scriptstyle \gamma_k\in\Gamma, 1 \leq k \leq N
       \end{matrix}} \tau\quad\quad\quad\text{(Flight time)}
\end{equation}
Such that we achieve:
\begin{enumerate}
    \item Correct synchronous trajectory:
        \[  \forall\,k\in[1,N] \quad \gamma_k(0) = S_k \ \wedge\ \gamma_k(\tau) = E_k\]
    \item Collision avoidance:
        \[  \forall\,k,s,\ k\neq s \quad \forall\,t\in[0,\tau]\ \norm{\gamma_k(t)-\gamma_s(t)} > \delta\]
    \item Constant speed:
        \[  \forall\,k \quad \forall\,t\in[0,\tau]\ \norm{\gamma_k'(t)} = V\]
    \item Bounded acceleration:
        \[  \forall\,k \quad \forall\,t\in[0,\tau]\ \norm{\gamma_k''(t)} \leq \frac{V^2}{\rho_\text{min}}\]
\end{enumerate}
The core starting point for solving this problem is to define the set of possible paths $\Gamma$.
Often this set is made of splines \cite{Yu_2017,Shanmugavel2007a,Cichella2014,Choe2016} as they provide
a continuous family with a convenient representation for integrating spatial and dynamic constraints
in a solver. But they heavily rely on a good initial guess and continuous optimization. This can be
tempered by the use of global heuristics like Particle Swarm Optimization (PSO), at a significant computational cost.

The method presented here is based on a finite set of solutions. 
It has the advantage of exploring significantly different paths,
which may then be used as a better starting point for locally optimized methods. 
We now present the elements making our set of solutions: the \emph{Dubins paths}.

\section{Dubins paths}
\label{sec:dubins}

Dubins showed in \cite{Dubins1957} that the shortest path between any two states of a vehicle described in \refeq{eq:dubins}
is made of straight lines and circle arcs of radius $\rho_\text{min}$. We denote the possible
components $S,L,R$, respectively for Straight, Left and Right [turn]. 
The optimal solution is among the following six possibilities: $RSR,LSL$ ; $RSL,LSR$ and $LRL,RLR$.
Since there is only a limited set of solutions, it is easy to find the best one by trying each of them%
.

\begin{figure}[b]
    \centering
    \includeinkscape[width=\columnwidth]{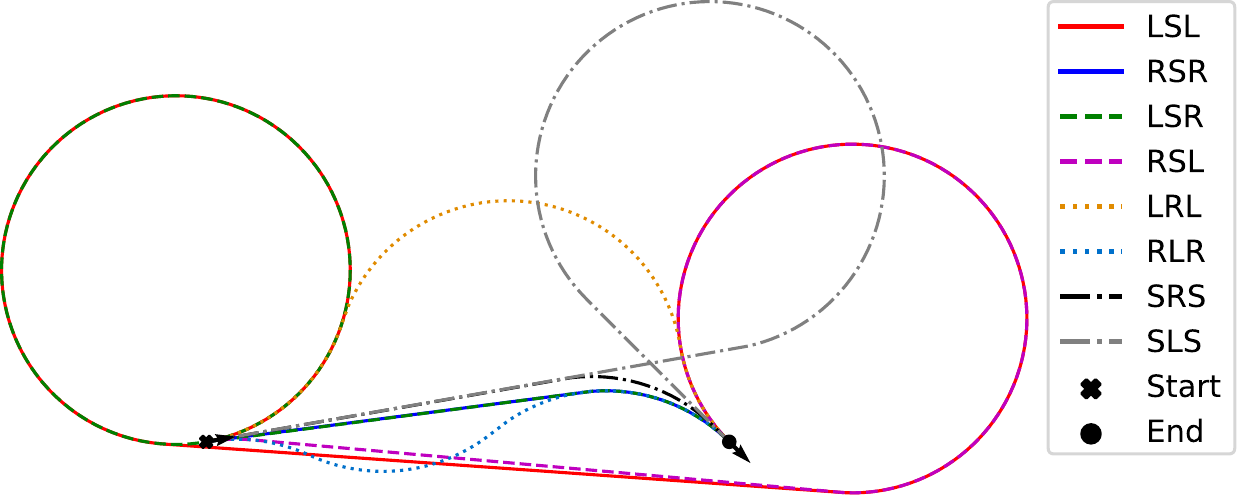}
    \caption{Examples of the six Dubins paths plus two extras (SRS, SLS)}
    \label{fig:dubins}
\end{figure}

We may still need more variety in the choice of paths in order to avoid conflicts.
We add the single-turn paths $SLS$ and $SRS$. They have a limited scope of application, but can still
extend the number of possibilities.
With the six basic Dubins paths and these two, we define the set of \emph{basic paths} $\Gamma_0$, which has cardinality 8.
For any type $\gamma\in\Gamma_0$, a path is exactly defined by starting and ending poses $S,E\in\setR^2\times\setS$ and a turn radius $\rho > 0$,
and we denote $\gamma[S,E,\rho]\in\mathcal{D}^2(\setR,\setR^2)$ the corresponding function, when it exists.
They are all showcased in \autoref{fig:dubins}.

By construction, all of these solutions are guaranteed to satisfy the Constant speed and Bounded acceleration constraints.
Two main questions remain: how to detect conflicts and how to tune the paths to achieve synchronicity ?

\subsection{Finding a fitting solution}

Dubins paths provide a direct solution to the individual path planning problem, but offer little tuning mechanisms
for achieving a specific path length. This question has been studied, and several methods devised:
artificially increase the turning radius \cite{Shanmugavel2005},
use circles with different radii \cite{Rao2024}, introduce straight extensions \cite{Chen2021} or deform straight parts \cite{Yao2017}.
We use two strategies: increasing the turn radius and adding straight segments at the start or end of the path.

More formally, fitting length $\ell > 0$ by tweaking radius $\rho$ of a path from $S = [x_S,y_S,\theta_S]$ to $E = [x_E,y_E,\theta_E]$
corresponds to the following optimization problem:
\begin{equation}
    \label{eq:rho-fit}
    \min_{\rho > \rho_\text{min}} (\text{len}(\gamma[S,E,\rho]) - \ell)^2
\end{equation}
Similarly, using length extensions distributed with ratio $r\in[0,1]$ along the initial direction $V_S = [\cos\theta_S,\sin\theta_S,0]$ 
and final one $V_E = [\cos\theta_E,\sin\theta_E,0]$  yields the problem:
\begin{equation}
    \label{eq:straight-fit}
    \min_{l \geq 0} (\text{len}(\gamma[S+l r V_S,E - l (1-r) V_E,\rho_\text{min}])+l - \ell)^2
\end{equation}

\begin{figure}
    \centering
    \includeinkscape[width=0.7\columnwidth]{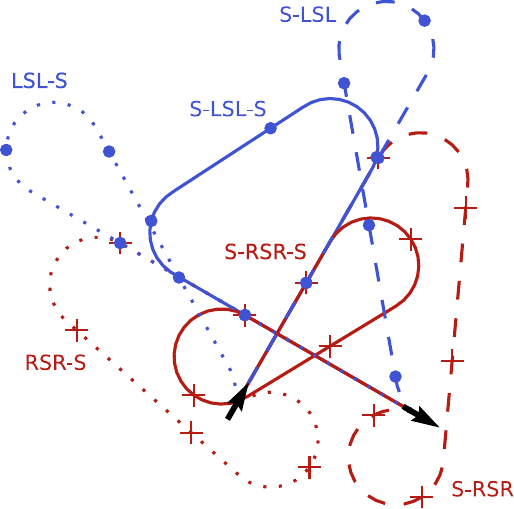}
    \caption{Start, end and both extended paths: $RSR$ and $LSL$}
    \label{fig:demo-extended-dubins}
\end{figure}

Note that the second method introduces new types of paths: Dubins extended by straight lines. The ratio parameter allows
many variations, but for the sake of simplicity we will only keep three: start extended ($r = 1$), noted S-\emph{AAA};
end extended ($r = 0$) \emph{AAA}-S and both extended ($r=\nicefrac{1}{2}$) S-\emph{AAA}-S, where in all cases \emph{AAA} represents
one of the three letters combination characterizing an element of $\Gamma_0$. A showcase of some Dubins extended paths
is given in \autoref{fig:demo-extended-dubins}. Let $\Gamma$ define the set containing 
the \emph{basic paths} $\Gamma_0$ length fitted using the aforementioned methods; its cardinal is $8+3\times 8$.

Solving problems \refeq{eq:rho-fit} and \refeq{eq:straight-fit} is manageable since they are unidimensional and
quite monotonous, in spite of some discontinuities (see \cite[Figure 4]{DubinsToBSplineSync}
for path length as a function of the turn radius). We use Brent's method \cite[Chapter 4]{Brent1973} to efficiently find solutions.

\subsection{Conflicts detection}

Dubins path are made of straight lines and circle arcs. Thus the question of detecting a conflict
between two paths of $\Gamma$ can be broken down into finding conflicts between their components. Three subcases emerge:
line-line, line-circle and circle-circle. Both 2D lines and circles can be conveniently represented using complex functions:
\begin{itemize}
    \item For a line, let $f(t) = v\cdot t + a$, with $v,a\in\setC$ and $\abs{v} = V$
    \item For a circle, $f(t) = C + r e^{i(\omega \cdot t + \phi)}$, 
    with $C\in\setC,\ \phi\in[0,2\pi],\ \omega\in\setR^+, r \geq \rho_\text{min}$ and $r \omega = V$
\end{itemize}

For two parametric functions $f_1,f_2\ \colon [0,\tau] \to \setC$ representing two trajectories,
we define two types of separations with respect to some minimal distance $\delta > 0$:
\begin{itemize}
    \item \emph{Temporal separation} : two points following prescribed trajectories must stay separated by at least $\delta$.
    \begin{equation}
        \label{eq:temp-sep}
        \min_{t\in[0,\tau]} \abs{f_1(t)-f_2(t)} > \delta 
    \end{equation}
    \item \emph{Spatial separation} : two shapes must be separated by at least $\delta$.
    \begin{equation}
        \label{eq:space-sep}
        \min_{(t_1,t_2)\in[0,\tau]^2} \abs{f_1(t_1)-f_2(t_2)} > \delta 
    \end{equation}
\end{itemize}
One can easily verify that if \refeq{eq:space-sep} is satisfied, then so is \refeq{eq:temp-sep}.

Practically, we are interested in \emph{temporal separation}. But, if it is easy to compute 
in the line-line case (it becomes a second degree polynomial equation),
this is not true for the line-circle and circle-circle cases. For those, even if it is a scalar optimization problem,
we require a global solver (the interval method \cite{Hansen1979} for instance). 
To limit the use of such computationally expensive procedures, we perform pre-filtering by computing the \emph{spatial separation}.
Indeed, since we are using only lines and circles, it is a straightforward geometry exercise to evaluate this separation.


We have presented in this section the set of paths $\Gamma$ used
for solving problem \refeq{eq:optim}. We have shown how to fit them to specific length
and how to detect conflicts between them. 
Let us now put this together into a fleet path planner.

\section{Fleet path planning}
\label{sec:fleetpp}

\subsection{Presentation of the algorithm}

Our procedure for planning a path for the fleet respecting the constraints given in
\refeq{eq:optim} works by combinatorial optimization. 
We consider two main factors regarding complexity:
the number of different possible paths $\abs{\Gamma}$ and the number of aircraft $N$. The former can be increased or
reduced by modifying the chosen individual path planning methods, and the latter expresses the scale and difficulty of the problem.

Three primitives are required to describe our method, whose underlying principles have been given in section \ref{sec:dubins}:
\begin{itemize}
    \item $\tau_k,\gamma_k \longleftarrow \mathrm{ShortestDubins}(S_k,E_k,V,\rho_\text{min})\ \colon$
    Given starting and ending states $S_k,E_k$, the aircraft speed $V$ and its minimal turn radius $\rho_\text{min}$,
    compute the shortest path $\gamma_k$ belonging to $\Gamma$ and give the travel time $\tau_k$.

    Its complexity scales linearly with $\abs{\Gamma}$, since each path must be considered to find the shortest.

    \item $L_k \longleftarrow\mathrm{FitDubins}(S_k,E_k,V,\rho_\text{min},\tau)\ \colon$
    Given starting and ending states $S_k,E_k$, the aircraft speed $V$, its minimal turn radius $\rho_\text{min}$
    and a travel time $\tau$, output a set of paths $L_k$ taken from $\Gamma$ with turning radii greater or equal than $\rho_\text{min}$
    and having a travel time of $\tau$. This set may be empty if there is no way to go from $S_k$ to $E_k$ given the
    constraints and possible paths.

    Its complexity also scales linearly with $\abs{\Gamma}$, since each path must be considered.

    \item $\mathrm{True}|\mathrm{False} \longleftarrow \mathrm{AreSeparated}(\delta,\gamma_1,\dots,\gamma_N)\ \colon$
    Given a minimal separation distance $\delta$ and a list of paths $\gamma_1,\dots,\gamma_N$, check
    if there is a loss of separation between any two pairs.

    We can safely assume that there is an underlying $\mathrm{IsPairSeparated}(\delta,\gamma_a,\gamma_b)$ subroutine
    checking for separation between two given paths, solving \emph{spatial} \refeq{eq:space-sep} then possibly \emph{temporal}
    \refeq{eq:temp-sep} separations. 
    $\mathrm{AreSeparated}$ requires $\frac{N(N-1)}{2} = \bigO{N^2}$ calls to the subroutine (since separation is a symmetric notion).
\end{itemize}

We can now describe the general principles, also summarized in pseudocode \autoref{alg:fleetpp}.
It starts by computing for each aircraft its best individual solution, ignoring collisions, by using the $\mathrm{ShortestDubins}$
primitive. We want synchronous arrivals, i.e. same travel time for everyone. 
A first good guess is the maximum of the minimum individual times, which we denote $\tau_\text{min}$. Using a user-defined parameter $R$,
the maximum possible time is set to $R\times \tau_\text{min}$. We use these two values to initialize a priority queue $Q$
containing the flight times considered.

The rest of \autoref{alg:fleetpp} executes the main loop until a stop condition is met:
For every untested flight time $\tau$ in $Q$, compute the possible paths $L_k$ for each aircraft $k\in[1,N]$ using $\mathrm{FitDubins}$,
and look through all combinations $(\gamma_1,\dots\gamma_N)\in L_1 \times \dots \times L_N$; if one satisfies $\mathrm{AreSeparated}$
with the desired separation $\delta$, break the loop, and otherwise move to the next value in $Q$.

If a solution was found, remove all elements in $Q$ with flight time greater than the one of the solution, since they cannot be optimal.
Then, $Q$ is populated by inserting $b$ values regularly spaced between every two successive $\tau_k,\tau_{k+1}\in Q$ (if $b=1$, this is equivalent
to successive dichotomies). To limit sampling to a reasonable resolution, a minimum width $w > 0$ is specified by the user, such that if
$\tau_{k+1} - \tau_k \leq w$, no samples are added between these two values.

The number of iterations is counted as the number of different flight times $\tau$ tested. Once above a given threshold, the program stops.
The second stop condition is maximum total execution time, and is evaluated before testing a new $\tau$. The last stop condition is progress;
the program stops if no new $\tau$ are added to $Q$ when re-sampling.

\begin{algorithm}
    \SetAlgoLined
    \DontPrintSemicolon
    \KwData{Initial and final states $(S_k,E_k)_{1\leq i\leq N}$; flight parameters $V,\rho_\text{min},\delta$; optimization parameters $R,w,b$}
    \KwResult{Duration $\tau$ and Dubins trajectories $(\gamma_k)_{1\leq i\leq N}$}
    \For(\tcp*[f]{Initial solutions}){$1\leq i\leq N$}
    {
        $\tau_k,\gamma_k \longleftarrow \mathrm{ShortestDubins}(S_k,E_k,V,\rho_\text{min})$\;
    }
    $\tau_\text{min} \longleftarrow \max_{1\leq i\leq N} (\tau_k)$\;
    \tcp{Ascending priority queue}
    $Q \longleftarrow [\tau_\text{min} ; R \cdot \tau_\text{min}]$\;
    Best $\longleftarrow \emptyset$\;
    \Repeat{\emph{Stop condition}}
    {
        \ForEach{untested $\tau\in Q$}
        {
            \tcp{Compute all solutions for a given length}
            \For{$1\leq i\leq N$}
            {
                $L_k \longleftarrow \mathrm{FitDubins}(S_k,E_k,V,\rho_\text{min},\tau)$\;
            }
            \tcp{If one $L_k$ is empty, this loop is skipped}
            \ForEach{Combination $(\gamma_1,\dots,\gamma_N)$ of $L_1 \times \dots \times L_N$}
            {
                \If{$\mathrm{AreSeparated}(\delta,\gamma_1,\dots,\gamma_N)$}
                {

                    Best $\longleftarrow \tau,(\gamma_1,\dots,\gamma_N)$\;
                    \Break twice\;
                }
            }
            Mark $\tau$ as tested\;
        }
        \If{\emph{Best} is not empty}
        {
            Remove from $Q$ all $\tau > $Best[0]\;
        }
        \For{Each successive pair $\tau_{i},\tau_{i+1} \in Q$}
        {
            $\Delta \tau \longleftarrow \tau_{i+1} - \tau_i$\;
            \If{$\Delta \tau > w$}
            {
                \lFor{$1 \leq j \leq b$}{Add $\tau_i + j\ \cfrac{\Delta\tau}{b+1}$ to $Q$} 
            }
        }
    }
    \Return Best\;
    \caption{Dubins based path planning algorithm}
    \label{alg:fleetpp}
\end{algorithm}

\subsection{Algorithm analysis}
\label{subsec:algo-complexity}

The different stop conditions ensure the program terminates in a timely manner, but does not guarantee
that a solution is found. Furthermore, there is no guarantee at large to find a solution, even with unlimited computations,
as it may not exist in the $\Gamma$ set chosen.
Still, we can evaluate the complexity of testing a flight duration $\tau$, which is the core of this algorithm.
Computing the set of possible paths for each aircraft contributes $N \cdot \bigO{\abs{\Gamma}}$.
Then, there are at most $\abs{\Gamma}$ paths per aircraft, meaning a total of $\abs{\Gamma}^N$ combinations.
Each combination has to be tested for conflicts, yielding a cost per iteration of $\bigO{N^2 \abs{\Gamma}^N}$.
This is \emph{polynomial} with respect to the number of different paths $\abs{\Gamma}$ and \emph{exponential} with respect to
the number of aircraft $N$.

There are ways to accelerate computations. Since conflicts involve aircraft two by two, it is possible to compute
for each of the $\frac{N(N-1)}{2}$ unordered aircraft pairs the $\abs{\Gamma}^2$ paths pairs, and check for conflict.
This provides $N^2 \abs{\Gamma}^2 $ binary coefficients $c_{a,b,g,h}$ denoting whether there is a conflict between
aircraft $a$ and $b$ using respectively paths $g$ and $h$. We can formulate a \emph{0-1 Integer Linear Program} (ILP)
with variables $x_{a,g}$ defining if aircraft $a$ uses path $g$. There are two types of constraints: avoiding conflicts
\[\forall a \neq b \in [1,N]\ \forall g,h\in\Gamma\ c_{a,b,g,h} \cdot (x_{a,g} + x_{b,h}) \leq 1\]
and choosing a unique path for each aircraft
\[\forall a\in[1,N]\ \sum_{g\in\Gamma} x_{a,g} = 1\]
This is a well-known NP-complete problem, and an efficient existing solver can be used, for instance \cite{Huangfu2018}.

Another important point is that this algorithm is easily parallelizable. Indeed, each aircraft $k$ can independently compute
its set of possible paths $L_k$, and each conflict check is independent of the others. Hence, if we consider a centralized
multiprocessor unit, we may use at most $\frac{\abs{\Gamma}^2 N(N-1)}{2}$ cores to compute the aforementioned $c_{a,b,g,h}$.

\subsection{Algorithm extensions}

The capabilities of \autoref{alg:fleetpp} can be extended by a clever use of its formulation.
The first one is that it is possible to take into account a constant uniform wind $W$. 
Since the flight time $\tau$ is set before planning the paths, it is possible to update the ending
point by integrating the wind deviation. This amounts to replacing $E_k$ by $E_k - W \tau$, thereby
anticipating the wind induced drift.
Note that since the same linear change of referential is applied to every aircraft,
conflict detection does not need to be changed.
Regarding the initial guess for $\tau_\text{min}$, there exist methods to compute it in presence of uniform wind,
for instance \cite{Wu2025}.

The second modification consists in incorporating temporal separation at arrival. Assume that aircraft $k\in[2,N]$
must reach their destination with a time difference $\Delta t_k$ with respect to the previous aircraft $k-1$.
This is implemented by replacing $\tau$ in the call of $\mathrm{FitDubins}$ by $\tau + \sum_{j=1}^{k-1} \Delta t_j$.
One use case of such modification is landing sequencing.

\section{Simulations}
\label{sec:results}

\begin{figure}
    \centering
    \includeinkscape[width=0.75\columnwidth]{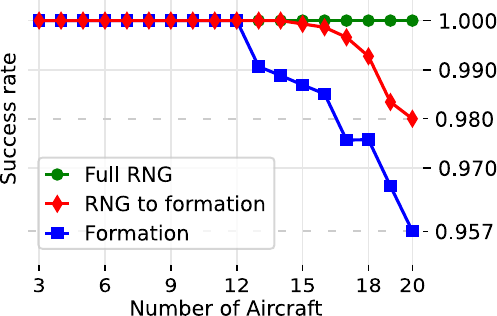}
    \caption{Success rate as function of aircraft number and type of test case}
    \label{fig:success-rate}
\end{figure}

The method is tested in simulation, based on three main types of scenarios: transition between formations,
random state to formation and random state to another random state, 
respectively referred as Formation, RNG to Formation and Full RNG in \autoref{fig:success-rate}.
Transitions between formations can be grouped into two subcategories: switch from a formation to a different one with same orientation,
or perform a turn with the same formation. A total of 22 formations were considered, ranging from circular and linear to rectangular and their variants.

A random initial state is generated by sampling points in an area, then using repulsion forces until minimal required separation is achieved, and
finally orientations are randomly picked. For a final state, it can be generated like the initial state, it can be uniformly shifted, or
it can be generated based on disks: for each initial point, its associated final point is sampled at a given distance from it, then
orientation is randomly picked. The repulsion forces are still applied to avoid trivially unfeasible cases. Note that this repulsion force may
push the points out of their initial area.

We based our simulations on a fixed-wing drone of \qty{1}{\meter} wingspan (equivalent to a Parrot Disco). The airspeed is set to \qty[per-mode = symbol]{15}{\meter\per\second},
minimal turn radius to \qty{40}{\meter} and minimal separation to \qty{80}{\meter}. 
When using random positions, initial sampling was done in a \qtyproduct[product-units = single]{1 x 1}{\kilo\meter}
square before applying repulsion until \qty{240}{\meter} of separation was achieved; this allows each aircraft to do a full circle
at minimal turn radius without loss of separation. 
When using positions in formation, the distance is set to \qty{120}{\meter}.
The distance between initial and final formations is of \qty{1}{\kilo\meter}. No wind was considered in these tests.
Regarding the algorithm parameters, we set $R = 3$, $b = 2$ and $w = \max(\qty{0.1}{\second} ,  R \times \tau_\text{min} \times 10^{-4})$.
The stop conditions values are 300 iterations maximum and a \qty{60}{\second} timeout.
The solver has been benchmarked on a 12th Gen Intel(R) Core(TM) i7-1255U, using 12 threads. 
An example of solutions produced by the solver is given in \autoref{fig:demo}, 
with two different wind conditions: 0 and \qty[per-mode = symbol]{10}{\meter\per\second}.

Tests have been done for fleets of sizes from 3 to 20 aircraft,
with 5230 cases each: 300 fully random, 530 formation transitions and 4400 random to formation. 
The proportion of successful cases per aircraft number and problem type is reported in \autoref{fig:success-rate}.
Success rate starts to decrease after 12 aircraft. 
Two main reasons can be given for failures. First, formation cases use tighter margins on endpoints and do not
optimize aircraft reordering between formations, instead using an arbitrary one which may be extra difficult. Second, randomly selected initial points
may hinder synchronicity, especially when the arrivals are tightly ordered. Nevertheless, for 12 aircraft and less,
no failures were reported, and the rate stays above 95\% when testing with 20 aircraft. 

\begin{figure}
    \centering
    \includeinkscape[width=\columnwidth]{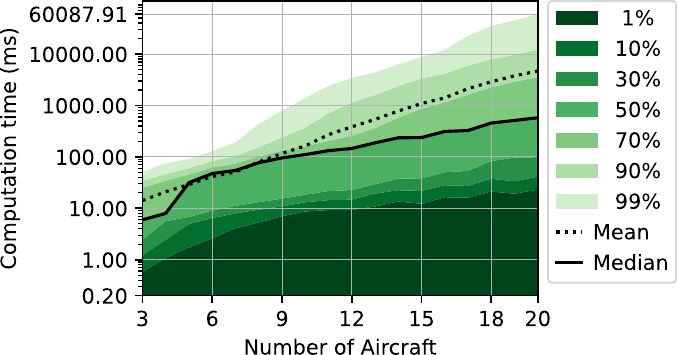}
    \caption{Computation time log distribution as function of aircraft number}
    \label{fig:compute-time}
\end{figure}

\begin{figure*}[t]
    \centering
    \includeinkscape[width=\textwidth]{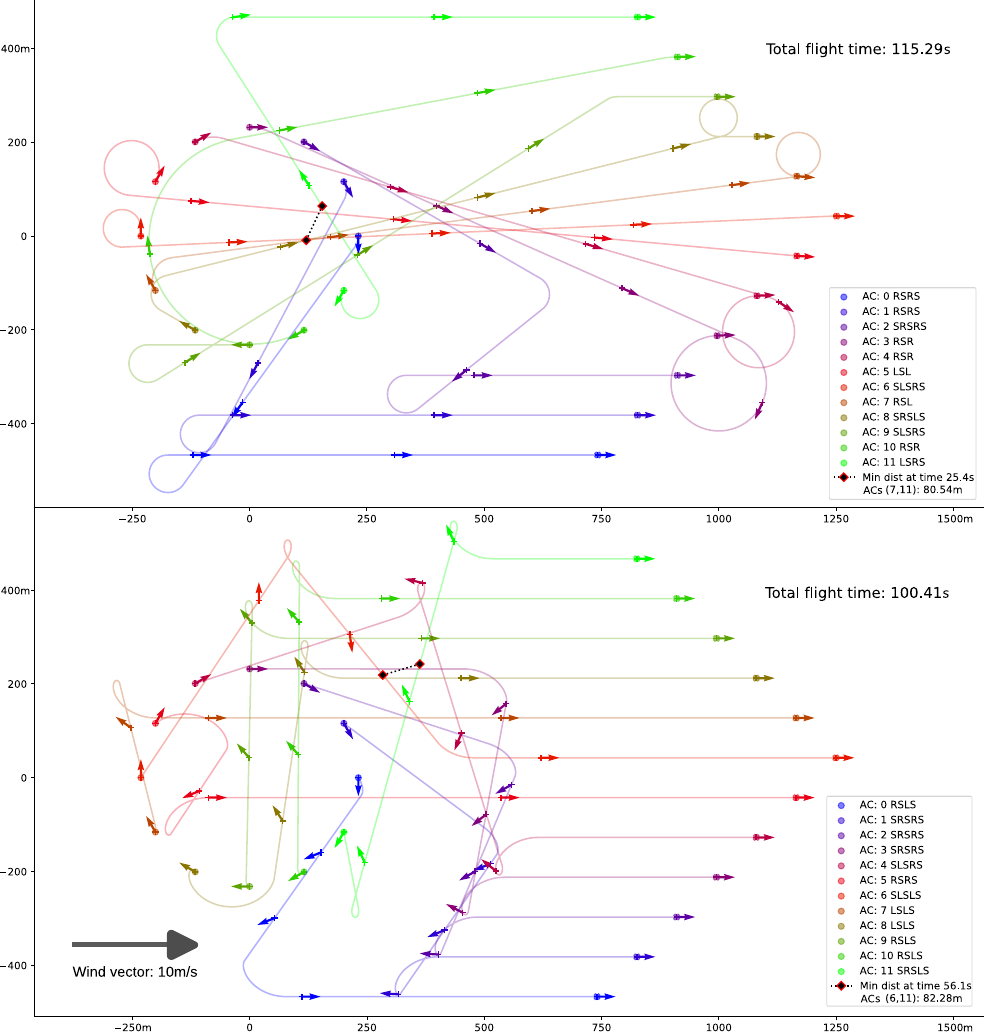}
    \caption{Demonstration of conflict free path planning without and with wind, from circular to chevron formations.}
    \label{fig:demo}
\end{figure*}

Computation times are reported in \autoref{fig:compute-time}. We measured the main thread time since parallelism has been exploited in our implementation.
For 9 aircraft and below, 99\% of cases are solved in less than \qty{1}{\second}. 
For 20 aircraft, 90\% of cases were each completed in less than \qty{13}{\second}. A detailed analysis reveals that the \qty{60}{\second} timeout is reached 
only for some cases with 17 aircraft or more. The number of iterations is never a limiting factor; the maximal value measured is 87 iterations.

The most difficult cases are the transitions between formations, then random to formations and finally fully random;
the computation times with 20 aircraft average to \qty{5065}{\milli\second}, \qty{3877}{\milli\second} and \qty{417}{\milli\second} respectively.
It can be explained by the lack of geometric separation between trajectories with formations, as it incurs the use of more time-consuming algorithms
for detecting potential conflicts.

\section{Conclusions and future works}

In this article we showed a method to solve synchronous wind-aware path planning for
a fleet of constant speed Dubins vehicles, and benchmarked it from 3 to 20 aircraft.
On average, even up to 20 aircraft, the computation time is \qty{8}{\second} for a success rate superior to 95\%.
The method enforces synchronicity by setting the travel time first, then creating a finite set of paths per aircraft satisfying the constraints,
and finally seeking a conflict-free combination, which is a solution to the fleet path planning problem.

This algorithm has several limitations, the first one being that there are no guarantees to find a solution. This may be mitigated
by adding possible paths and reducing margins, but these come with computational and safety costs.
Secondly, computational cost scales exponentially with the number of aircraft, limiting scalability.
This method was not made to replace swarm-like tactics used to handle large amounts of UAVs, but the
time to find a solution can be an obstacle depending on the available hardware and expected reactivity.
Lastly, the paths presented here are based on Dubins model. If it simplifies greatly path planning, it does not
match an actual aircraft model, since lateral acceleration is discontinuous. Given how important separation
between aircraft is, we wish that path following is as precise as possible. 

There are multiple ways of improving this method. One could integrate vertical movement \cite{Owen2014},
or incorporate static obstacle avoidance \cite{DubinsObstacle}.
Computation time could be improved by studying how conflicts between paths evolve as they are continuously deformed
to fit a different length, by using homotopy based methods for instance \cite{VazquezLeal2013}.
Such insight could limit the need to recompute the ILP coefficients defined in \autoref{subsec:algo-complexity}.
Lastly, it is possible to locally improve the solution by smoothing without violating the satisfied constraints,
using B-Splines \cite{DubinsToBSplineSync} or Pythagorean Hodograph splines \cite{Shanmugavel2007a,Choe2016}.

Our C++ implementation of the solver can be found at \url{https://github.com/enacuavlab/DubinsFleetPlanner}.






\printbibliography

\end{document}